\documentstyle[a4,12pt,amsfonts,leqno]{article}

\newtheorem{Theorem}{Theorem}[section]

\newtheorem{Lemma}[Theorem]{Lemma}
\newtheorem{Corollary}[Theorem]{Corollary}
\newtheorem{Remark}[Theorem]{Remark}
\newtheorem{Example}[Theorem]{Example}

\newtheorem{Question}[Theorem]{Question}

\newcounter{ctr}[section]

\newcommand{\bangou}{\addtocounter{ctr}{1}}

\newcommand{\rmn}[1]{\romannumeral#1}
\def\RMN#1{\uppercase\expandafter{\romannumeral#1}}

\newcommand{\ya}[1]{\mathrel{\overrightarrow{\hphantom{#1}}}}

\newcommand{\qed}{{\unskip\nobreak\hfil\penalty50\quad\null\nobreak\hfil{\bf
q.e.d.}\parfillskip0pt\finalhyphendemerits0\par\medskip}}
\newcommand{\proof}{\noindent{\it Proof.} \ }

\newcommand{\rank}{\mathop{\rm rank}\nolimits}

\newcommand{\spec}{\mathop{\rm Spec}\nolimits}

\newcommand{\codim}{\mathop{\rm codim}\nolimits}

\newcommand{\g}{\mathop{\rm G}\nolimits_0}

\newcommand{\chow}[1]{\mathop{\rm A}\nolimits_{#1}}
\newcommand{\subq}{_{\Bbb Q}}

\newcommand{\p}{{\frak p}}
\newcommand{\q}{{\frak q}}
\newcommand{\m}{{\frak m}}

\begin{document}
\title{The singular Riemann-Roch theorem and Hilbert-Kunz functions}
\author{Kazuhiko Kurano\thanks{The author is
supported by a Grant-in-Aid for scientific Research Japan.} (Meiji University)}
\date{}
\maketitle

\begin{abstract}
In the paper, via the singular Riemann-Roch theorem, it is proved that the class of
the $e$-th Frobenius power ${}^eA$ can be described using the class of the canonical module 
$\omega_A$ for a normal local ring $A$ of positive characteristic.
As a corollary, we prove that the coefficient $\beta(I,M)$ of
the second term of the 
Hilbert-Kunz function $\ell_A(M/I^{[p^e]}M)$ of $e$ vanishes if $A$ is a ${\Bbb Q}$-Gorenstein ring and $M$ is a finitely generated $A$-module 
of finite projective dimension.

For a normal algebraic variety $X$ over a perfect field of positive characteristic,
it is proved that the first Chern class of the $e$-th Frobenius power $F^e_*{\cal O}_X$
can be described using the canonical divisor $K_X$.
\end{abstract}

\section{Introduction}
Let $(A, \m)$ be a $d$-dimensional Noetherian local ring of characteristic $p$,
where $p$ is a prime integer.
Here, $\m$ is the unique maximal primary ideal of $A$.
For an $\m$-primary ideal $I$ and a positive integer $e$, we set
\[
I^{[p^e]} = (a^{p^e} \mid a \in I)A .
\]
It is easy to see that $I^{[p^e]}$ is an $\m$-primary ideal of $A$.
For a finitely generated $A$-module $M$, 
the function $\ell_A(M/I^{[p^e]}M)$ of $e$ is called the {\em Hilbert-Kunz function}
of $M$ with respect to $I$, where $\ell_A( \ )$ stands for the length of the given
$A$-module.
It is known that 
\[
\lim_{e \rightarrow \infty} \frac{\ell_A(M/I^{[p^e]}M)}{p^{de}}
\]
exists~\cite{Mo}, and this limit is called the {\em Hilbert-Kunz multiplicity},
which is denoted by $e_{HK}(I,M)$.
Several properties of $e_{HK}(I,M)$ have been studied by many authors (Monsky, Watanabe, Yoshida, Huneke, Enescu, etc.).

Recently Huneke, McDermott and Monsky (Theorem~1, Corollary~1.10 and Theorem~1.11 in \cite{HMM}) proved the following exciting theorem:

\bangou
\begin{Theorem}[Huneke, McDermott and Monsky]\label{HMM}
Let $(A, \m)$ be a $d$-dimensional excellent normal local ring 
of characteristic $p$, where $p$ is a prime integer.
Assume that the residue class field of $A$ is perfect.

Let $I$ be an $\m$-primary ideal of $A$ and $M$ be a finitely generated $A$-module.
\begin{enumerate}
\item
There exists a real number $\beta(I,M)$ that satisfies the following equation:\footnote{
Let $f(e)$ and $g(e)$ be functions of $e$. 
We denote $f(e) = O(g(e))$ if there exists a real number $K$ that satisfies
$| f(e) | < K g(e) $ for all $e \gg 0$.
}
\[
\ell_A(M/I^{[p^e]}M) = e_{HK}(I,M) \cdot p^{de} + \beta(I,M) \cdot p^{(d-1)e}
+ O(p^{(d-2)e}) 
\]
\item
Assume that $A$ is F-finite.\footnote{
We say that $A$ is F-finite if the Frobenius map $F : A \rightarrow A = {}^1A$ 
is module-finite.
We sometimes denote the $e$-th iteration of $F$ by 
$F^e : A \rightarrow A = {}^eA$.}
Then, there exists a ${\Bbb Q}$-homomorphism $\tau_I : {\rm Cl}(A)\subq \longrightarrow 
{\Bbb R}$ that satisfies
\[
\beta(I,M) = \tau_I
\left(
{\rm cl}(M) - \frac{\rank_AM}{p^d - p^{d-1}} {\rm cl}({}^1A)
\right) 
\]
for any finitely generated torsion-free $A$-module $M$.
In particular, we have
\[
\beta(I,A) = - \frac{1}{p^d - p^{d-1}} \tau_I
\left(
 {\rm cl}({}^1A)
\right) .
\]
\end{enumerate}
\end{Theorem}

We denote by ${\Bbb Q}$ (resp.\ ${\Bbb R}$) the field of rational numbers
(resp.\ real numbers).
For an abelian group $N$, 
$N\subq$ stands for $N \otimes_{\Bbb Z} {\Bbb Q}$.

The map ${\rm cl} : \g(A) \longrightarrow {\rm Cl}(A)$ is defined by Bourbaki~\cite{B}
and sometimes called the {\em determinant map} (see Remark~\ref{rem} below).

It is natural to ask the following questions:

\bangou
\begin{Question}
\begin{rm}
\begin{enumerate}
\item
When does ${\rm cl}({}^1A)$ vanish?
\item
How do the ${\rm cl}({}^eA)$'s behave?
\end{enumerate}
\end{rm}
\end{Question}

Using the singular Riemann-Roch formula, we obtain the following theorem:

\bangou
\begin{Theorem}\label{main}
Let $(A, \m)$ be a $d$-dimensional 
Noetherian normal local ring of characteristic $p$,
where $p$ is a prime integer.
Assume the following three conditions;
(1) \ $A$ is a homomorphic image of a 
regular local ring, (2) the residue class field of $A$ is perfect, 
and (3) $A$ is F-finite.

Then, for each integer $e > 0$, we have
\[
{\rm cl}({}^eA) = \frac{p^{de} - p^{(d-1)e}}{2}{\rm cl}(\omega_A)
\]
in ${\rm Cl}(A)_{\Bbb Q}$, where $\omega_A$ is the canonical module of $A$.
\end{Theorem}

The following corollary is an immediate consequence of 
Theorem~\ref{HMM} and Theorem~\ref{main}.
A Noetherian normal local domain $A$ is called ${\Bbb Q}$-{\em Gorenstein}
if ${\rm cl}(\omega_A)$ is a torsion element in ${\rm Cl}(A)$.

\bangou
\begin{Corollary}\label{cor}
Under the same assumptions as in Theorem~\ref{main},
if $A$ is a ${\Bbb Q}$-Gorenstein ring,
then $\beta(I,A) = 0$ for any maximal primary ideal $I$.
\end{Corollary}

\bangou
\begin{Remark}
\begin{rm}
If $A$ is a ${\Bbb Q}$-Gorenstein ring, then we have 
$\beta(I,M) = \tau_I({\rm cl}(M))$ by Theorem~\ref{HMM} (2).

Furthermore, assume that $M$ is a finitely generated $A$-module of
finite projective dimension.
Then, we have ${\rm cl}(M) = {\rm rank}_AM \cdot {\rm cl}(A) = 0$.
Therefore, in this case, $\beta(I,M)$ is equal to $0$.
\end{rm}
\end{Remark}

The following is an analogue of Theorem~\ref{main} 
for normal algebraic varieties.

\bangou
\begin{Theorem}\label{analogue}
Let $k$ be a perfect field of characteristic $p$,
where $p$ is a prime integer.
Let $X$ be a normal algebraic variety over $k$ of dimension $d$.
Let $F : X \rightarrow X$ be the absolute Frobenius map.\footnote{
Remark that, under the assumption, $F$ is a finite morphism.}

Then, we have 
\[
c_1(F_*^e{\cal O}_X) = \frac{p^{de} - p^{(d-1)e}}{2} K_X
\]
in $\chow{d-1}(X)_{\Bbb Q}$,
where $c_1( \ )$ is the first Chern class\footnote{
Set $U = X \setminus {\rm Sing}(X)$.
Since $\codim_X{\rm Sing}(X) \geq 2$, 
the restriction $\chow{d-1}(X) \rightarrow \chow{d-1}(U)$ is an isomorphism.
Here, remark that $(F_*^e{\cal O}_X)|_U = (F|_U)^e_*{\cal O}_U$ is 
a locally free sheaf on $U$.
Thus, $c_1(F_*^e{\cal O}_X)$ is defined as the first Chern class
$c_1((F_*^e{\cal O}_X)|_U) \in \chow{d-1}(U) = \chow{d-1}(X)$.
}
and $K_X$ is the canonical divisor of $X$.
\end{Theorem}

Here, $A_{d-1}(X)$ is the Chow group of $X$ consisting of cycles of dimension $d-1$.
We refer the reader to \cite{F} for Chow groups.
If $A$ (in Theorem~\ref{main}) is a local ring at a closed point 
of a normal algebraic variety over a perfect
field of positive characteristic, then Theorem~\ref{main} follows from
Theorem~\ref{analogue}.

We give a proof of Theorem~\ref{main} and Theorem~\ref{analogue} 
in the next section.

\section{A proof of Theorem~\ref{main} and Theorem~\ref{analogue}}

Before proving Theorem~\ref{main}, we recall basic properties on 
the determinant map.

\bangou
\begin{Remark}\label{rem}
\begin{rm}
The map ${\rm cl}$ in Theorem~\ref{main} is called the {\em determinant map}
which is defined by Bourbaki~\cite{B}.
Here, recall basic properties on ${\rm cl}$
which are used later.

Let $R$ be a Noetherian normal domain.
The group of isomorphism classes of reflexive $R$-modules 
of rank $1$ is called 
the {\em divisor class group} of $R$, and denoted by ${\rm Cl}(R)$.
Let $\g(R)$ be the Grothendieck group of finitely generated $R$-modules.
For an $R$-module $M$, we denote by $[M]$ the element in $\g(R)$
corresponding to the isomorphism class which $M$ belongs to.
Then, there exists the unique map 
\[
{\rm cl} : \g(R) \longrightarrow {\rm Cl}(R)
\]
that satisfies the following two conditions:
\begin{itemize}
\item[(1)]
If $M$ is a reflexive module of rank $1$, then ${\rm cl}([M])$ is just the isomorphism class which $M$ belongs to.
\item[(2)]
Let $M$ be a finitely generated $R$-module.
If the height of the annihilator of $M$ is greater than $1$,
then ${\rm cl}([M]) = 0$.
\end{itemize}
For an $R$-module $M$, we denote ${\rm cl}([M])$ simply by ${\rm cl}(M)$ 
as usual.
\end{rm}
\end{Remark}
%
%

Let $(A,\m)$ be a Noetherian local ring that satisfies the assumption
in Theorem~\ref{main}.
It is enough to prove Theorem~\ref{main} for complete local rings.
Therefore, in the rest of this section, 
we assume that $A$ is a $d$-dimensional local normal domain which is a homomorphic image of 
a formal power series ring $S$ over a perfect field $k$ 
of positive characteristic
unless otherwise specified.
By the singular Riemann-Roch theorem (cf., Chapter~18 and 20 in \cite{F}),
we have an isomorphism 
\[
\tau_{\spec(A)/\spec(S)} : \g(A)\subq \longrightarrow \chow{*}(A)\subq
\]
of ${\Bbb Q}$-vector spaces.
Here, remark that the Riemann-Roch map $\tau_{\spec(A)/\spec(S)}$ is determined 
not only by $\spec(A)$ but also by the regular base scheme 
$\spec(S)$ as in 20.1 in \cite{F}.
Let
\[
p_i : \chow{*}(A)\subq \longrightarrow \chow{i}(A)\subq 
\]
be the projection for $i = 0, 1, \ldots, d$.
We set 
\[
\tau_{i} = p_i \circ \tau_{\spec(A)/\spec(S)} : 
\g(A)\subq \longrightarrow \chow{i}(A)\subq .
\]
For a prime ideal $\p$ of $A$, $[\spec(A/\p)]$ stands for the element
in $\chow{*}(A)$ corresponding to the closed subscheme
$\spec(A/\p)$ of $\spec(A)$.

\bangou
\begin{Lemma}\label{BF}
Let $A$ be a $d$-dimensional normal local ring which is 
a homomorphic image of a regular local ring.
\begin{itemize}
\item[(\rmn{1})]
There exists a natural isomorphism $\chow{d-1}(A) = {\rm Cl}(A)$ by
identifying $[\spec(A/\p)]$ with ${\rm cl}(\p)$
for any prime ideal $\p$ of height $1$.
Then, for any prime ideal $\q \neq 0$, 
$\tau_{d-1}([A/\q])$ is equal to $- {\rm cl}(A/\q)$.
\item[(\rmn{2})]
We have the equality
\[
\tau_{d-1}([A]) = \frac{1}{2} {\rm cl}(\omega_A)
\]
in $\chow{d-1}(A)\subq = {\rm Cl}(A)\subq$. 
\item[(\rmn{3})]
Furthermore, assume that $A$ is a homomorphic image of a formal power series ring $S$ over a perfect field of positive characteristic.
Then, for each $e > 0$ and $i = 0, 1, \ldots, d$, the equality
\[
\tau_{i}([{}^eA]) = p^{ie}\tau_{i}([A]) 
\]
is satisfied.
\end{itemize}
\end{Lemma}

\proof
First we prove (\rmn{1}).
It is well-known that there exists an isomorphism 
$\chow{d-1}(A) \rightarrow {\rm Cl}(A)$ by 
$[\spec(A/\p)] \mapsto {\rm cl}(\p)$ (cf., Bourbaki~\cite{B}).
Suppose that $\p$ is a prime ideal of height $1$.
By the exact sequence
\[
0 \rightarrow \p \rightarrow A \rightarrow A/\p \rightarrow 0 ,
\]
we have 
\[
{\rm cl}(\p) = {\rm cl}(A) - {\rm cl}(A/\p) = - {\rm cl}(A/\p) .
\]
On the other hand, by the top-term property (Theorem~18.3 (5) in \cite{F}),
we have $\tau_{d-1}([A/\p]) = [\spec(A/\p)]$.
Thus, we obtain
\[
\tau_{d-1}([A/\p]) = [\spec(A/\p)] = {\rm cl}(\p) = - {\rm cl}(A/\p) .
\]
Let $\q$ be a prime ideal of height at least $2$.
By the top-term property, we have $\tau_{d-1}([A/\q]) = 0$.
In this case, we also have ${\rm cl}(A/\q) = 0$ by Remark~\ref{rem} (2).
The proof of (\rmn{1}) is completed.

We refer the reader to Lemma~3.5 of \cite{K9} for a proof of (\rmn{2}).

Now we start to prove (\rmn{3}).
Consider the following commutative diagrams:
\[
\begin{array}{ccc}
{}^eS & \longrightarrow & {}^eA \\
\scriptstyle{F^e}{\displaystyle \uparrow}\phantom{\scriptstyle{F^e}} & & 
\phantom{\scriptstyle{F^e}}{\displaystyle \uparrow}\scriptstyle{F^e} \\
S & \longrightarrow & A 
\end{array}
\ \ \ \ \ \ \ 
\begin{array}{ccccc}
S^{p^e} & \longrightarrow & S & \longrightarrow & A \\
\parallel & & & & \parallel \\
S & \longrightarrow & A & \stackrel{F^e}{\longrightarrow} & {}^eA
\end{array}
\]
The lefthand diagram above and the covariance with the proper map
$F^e : \spec({}^eA) \rightarrow \spec(A)$ (Theorem~18.3 (1) in \cite{F})
imply that the bottom half of the following diagram commutes.
The righthand diagram above implies that 
the top half of the following diagram commutes.
\bangou
\begin{equation}\label{kakan}
\begin{array}{ccc}
\g(A)\subq & \stackrel{\tau_{\spec(A)/\spec(S^{p^e})}}{\ya{aaaaaaaaaa}} &
\chow{*}(A)\subq \\
\parallel & & \parallel \\
\g({}^eA)\subq & \stackrel{\tau_{\spec({}^eA)/\spec(S)}}{\ya{aaaaaaaaaa}} &
\chow{*}({}^eA)\subq \\
\scriptstyle{F^e_*} {\displaystyle \downarrow} \phantom{\scriptstyle{F^e_*}}
& & 
\scriptstyle{F^e_*} {\displaystyle \downarrow} \phantom{\scriptstyle{F^e_*}}
\\
\g(A)\subq & 
\stackrel{\tau_{\spec(A)/\spec(S)}}{\ya{aaaaaaaaaa}} & \chow{*}(A)\subq
\end{array} 
\end{equation}
Here, $S^{p^e} = \{ x^{p^e} \mid x \in S \} \subset S$.
Remark that $S^{p^e}$ is a regular local ring and $S$ is a finite module 
over $S^{p^e}$.
Therefore, $\tau_{\spec(A)/\spec(S^{p^e})}$ and 
$\tau_{\spec(S)/\spec(S^{p^e})}$ can be defined (cf.\ Chapter~18 and 20
in Fulton~\cite{F}).

Here, we shall prove 
\bangou
\begin{equation}\label{claim}
\tau_{\spec(A)/\spec(S^{p^e})} = \tau_{\spec(A)/\spec(S)} 
\end{equation}
for any $e > 0$.

Since $S$ is a regular local ring, we have $\g(S)\subq = {\Bbb Q}[S]$
and $\chow{*}(S)\subq = {\Bbb Q}[\spec(S)]$
since $\g(S)\subq$ is isomorphic to $\chow{*}(S)\subq$ by the singular 
Riemann-Roch theorem.
By the top term property (Theorem~18.3 (5) in \cite{F}), we have
\bangou
\begin{equation}\label{topterm}
\tau_{\spec(S)/\spec(S^{p^e})}([S]) = [\spec(S)] .
\end{equation}

Let $M$ be a finitely generated $A$-module and
${\Bbb F}.$ be a finite $S$-free resolution of $M$.
By definition of $\tau_{\spec(A)/\spec(S)}$ (18.3 in \cite{F}),
we have
\bangou
\begin{equation}\label{A'ÆS}
\tau_{\spec(A)/\spec(S)}([M]) = 
{\rm ch}^{\spec(S)}_{\spec(A)}({\Bbb F}.) \cap [\spec(S)]
\end{equation}
in $\chow{*}(A)\subq$,
where 
\[
{\rm ch}^{\spec(S)}_{\spec(A)}({\Bbb F}.) :
\chow{*}(S)\subq \longrightarrow \chow{*}(A)\subq
\]
is the localized Chern character of the complex ${\Bbb F}.$ (18.1 in \cite{F}).
Therefore, we have
\begin{eqnarray*}
\tau_{\spec(A)/\spec(S^{p^e})}([M]) & = & 
\tau_{\spec(A)/\spec(S^{p^e})} \left(
\sum_i (-1)^i [H_i({\Bbb F}.)]
\right) \\  
& = & {\rm ch}^{\spec(S)}_{\spec(A)}({\Bbb F}.) \cap 
\tau_{\spec(S)/\spec(S^{p^e})}([S]) \\
& = & {\rm ch}^{\spec(S)}_{\spec(A)}({\Bbb F}.) \cap [\spec(S)] \\
& = & \tau_{\spec(A)/\spec(S)}([M]) ,
\end{eqnarray*}
where the second equality follows from the local Riemann-Roch formula 
(Example~18.3.12 in \cite{F}),
the third from (\ref{topterm}) and the fourth from (\ref{A'ÆS}).
Thus, (\ref{claim}) has been proved.

We denote the composite maps of the vertical arrows of 
the diagram (\ref{kakan})
by $F^e_* : \g(A)\subq \rightarrow \g(A)\subq$
and
$F^e_* : \chow{*}(A)\subq \rightarrow \chow{*}(A)\subq$, respectively.
Then, by the definition of $F_*$ (cf., Chapter~1 in \cite{F}), it is easy to see that 
the restriction $F^e_*|_{\chow{i}(A)\subq}$ is just the multiplication by $p^{ie}$
for $i = 0, 1, \ldots, d$ and any $e > 0$.
By (\ref{kakan}) and (\ref{claim}), we have the following commutative diagram:
\bangou
\begin{equation}\label{comm}
\begin{array}{ccc}
\g(A)\subq & \stackrel{\tau_{\spec(A)/\spec(S)}}{\ya{aaaaaaaaaa}} &
\chow{*}(A)\subq \\
\scriptstyle{F^e_*} {\displaystyle \downarrow} \phantom{\scriptstyle{F^e_*}}
& & 
\scriptstyle{F^e_*} {\displaystyle \downarrow} \phantom{\scriptstyle{F^e_*}}
\\
\g(A)\subq & 
\stackrel{\tau_{\spec(A)/\spec(S)}}{\ya{aaaaaaaaaa}} & \chow{*}(A)\subq
\end{array} 
\end{equation}
Thus, for an $A$-module $M$,
we obtain 
\bangou
\begin{equation}\label{rmn3}
\tau_iF^e_*([M]) = p^{ie} \tau_i([M])
\end{equation}
in $\chow{i}(A)\subq$ for each $i$ and $e$.
Since $F^e_*([A]) = [{}^eA]$, (\rmn{3}) has been proved.
\qed

Before proving Theorem~\ref{main}, we prove the following lemma:

\bangou
\begin{Lemma}\label{lem}
Let $(A,\m)$ be a $d$-dimensional 
normal local ring that is a homomorphic image of a regular 
local ring.
Then, for a finitely generated $A$-module $M$, we have
\[
\tau_{d-1}([M]) = -{\rm cl}(M) + \frac{\rank_AM}{2} {\rm cl}(\omega_A) 
\]
in ${\rm Cl}(A)\subq$.
\end{Lemma}

\proof
Set $r = \rank_AM$.
Then we have an exact sequence
\[
0 \rightarrow A^r \rightarrow M \rightarrow T \rightarrow 0 ,
\]
where $T$ is a torsion module.
By this exact sequence, we obtain
\[
{\rm cl}(M) = r \cdot {\rm cl}(A) + {\rm cl}(T) = {\rm cl}(T) .
\]
On the other hand, by (\rmn{2}) in Lemma~\ref{BF},
we obtain
\[
\tau_{d-1}([M]) = r \cdot \tau_{d-1}([A]) + \tau_{d-1}([T]) = 
\frac{r}{2}{\rm cl}(\omega_A) + \tau_{d-1}([T]) .
\]

Therefore, it is enough to show $\tau_{d-1}([T]) = - {\rm cl}(T)$ for any torsion module $T$.

We may assume that $T = A/\q$, where $\q \neq 0$ is a prime ideal of $A$.
Then, by (\rmn{1}) in Lemma~\ref{BF}, we have 
\[
\tau_{d-1}([A/\q]) = - {\rm cl}(A/\q) 
\]
as required.
\qed

Now we start to prove Theorem~\ref{main}.

By (\rmn{3}) and (\rmn{2}) in Lemma~\ref{BF}, we obtain
\[
\tau_{d-1}([{}^eA]) = p^{(d-1)e}\tau_{d-1}([A]) = \frac{p^{(d-1)e}}{2} {\rm cl}(\omega_A) .
\]
By Lemma~\ref{lem}, we have
\[
\tau_{d-1}([{}^eA]) = - {\rm cl}({}^eA) + \frac{\rank_A{}^eA}{2} {\rm cl}(\omega_A)
\]
in ${\rm Cl}(A)\subq$.
Since $\rank_A{}^eA = p^{de}$,
we obtain
\[
{\rm cl}({}^eA) = \frac{p^{de} - p^{(d-1)e}}{2} {\rm cl}(\omega_A)
\]
in ${\rm Cl}(A)\subq$.
\qed

\bangou
\begin{Remark}
\begin{rm}
By Lemma~\ref{lem} and Theorem~\ref{main}, we have
\[
\tau_{d-1}([M]) = 
- {\rm cl}(M) + \frac{\rank_AM}{2} {\rm cl}(\omega_A)
= - {\rm cl}(M) + \frac{\rank_AM}{p^d-p^{d-1}} {\rm cl}({}^1A) .
\]
Therefore, by Theorem~\ref{HMM}, we have 
\[
\beta(I,M) = -\tau_I(\tau_{d-1}([M])) \ \ 
\mbox{and} \ \ 
\beta(I,A) = -\frac{1}{2}\tau_I({\rm cl}(\omega_A)) 
\]
for any torsion-free $A$-module $M$.
\end{rm}
\end{Remark}

In the rest of this section, we shall give an outline 
of a proof of Theorem~\ref{analogue}.
Let $X$ be an algebraic variety that satisfies the assumptions in Theorem~\ref{analogue}.
Removing singularities of $X$, we may assume that $X$ is a smooth algebraic variety
over a perfect field $k$ of characteristic $p > 0$.\footnote{
Set $U = X \setminus {\rm Sing}(X)$.
Since $\codim_X{\rm Sing}(X) \geq 2$, 
the restriction $\chow{d-1}(X) \rightarrow \chow{d-1}(U)$ is an isomorphism.
On the other hand, we have $(F_*^e{\cal O}_X)|_U = (F|_U)^e_*{\cal O}_U$
and $K_X|_U = K_U$.
Therefore, we have only to show 
$c_1((F|_U)_*^e{\cal O}_U) = \frac{p^{de} - p^{(d-1)e}}{2} K_U$.

In the case of Theorem~\ref{main},
the proof does not become easier even if we remove singularities of 
$\spec(A)$.
The reason is that $\spec(A) \setminus {\rm Sing}(A)$ is not smooth over 
the base regular scheme $\spec(S)$.
}
Applying much the same method as in the proof of the commutativity of the diagram 
(\ref{comm}), one can prove that the diagram
\[
\begin{array}{ccc}
\g(X)\subq & \stackrel{\tau_{X/\spec(k)}}{\longrightarrow} 
& 
\chow{*}(X)\subq
 \\
\scriptstyle{F^e_*} {\displaystyle \downarrow} \phantom{\scriptstyle{F^e_*}} & 
&
\scriptstyle{F^e_*} {\displaystyle \downarrow} \phantom{\scriptstyle{F^e_*}} \\
\g(X)\subq & 
\stackrel{\tau_{X/\spec(k)}}{\longrightarrow} 
& \chow{*}(X)\subq \\
\end{array} 
\]
is also commutative.

Set 
\[
\tau_{X/\spec(k)}([{\cal O}_X]) = t_d + t_{d-1} + \cdots + t_0
\]
where $t_i \in \chow{i}(X)\subq$ for $i = 0, 1, \cdots, d$.
By the commutative diagram above, we have
\bangou
\begin{equation}\label{eq1}
\tau_{X/\spec(k)}([F^e_*{\cal O}_X]) = 
p^{de}t_d + p^{(d-1)e}t_{d-1} + \cdots + p^0t_0 .
\end{equation}
On the other hand, by Theorem~18.3 (2) in \cite{F}, we have
\bangou
\begin{eqnarray}
\nonumber
\tau_{X/\spec(k)}([F^e_*{\cal O}_X]) &  =  & {\rm ch}(F^e_*{\cal O}_X) 
\cap \tau_{X/\spec(k)}([{\cal O}_X]) \\
\label{eq2}
 & = & \left( p^{de} + c_1 + \frac{1}{2}(c_1^2 - 2c_2) + \cdots   \right)
 \cap \left( t_d + t_{d-1} + t_{d-2} + \cdots \right)
\end{eqnarray}
where $c_i$ stands for $c_i(F^e_*{\cal O}_X)$ for $i = 1, 2, \ldots$
(cf. Example~3.2.3 in \cite{F}).
Here, remark that $F^e_*{\cal O}_X$ is a locally free sheaf on $X$ 
since $X$ is a non-singular variety.
Comparing (\ref{eq1}) with (\ref{eq2}), we have the following equalities:
\bangou
\begin{eqnarray}
\nonumber
p^{de}t_d & = & p^{de}t_d \\
\label{eq3}
p^{(d-1)e}t_{d-1} & = & c_1t_d + p^{de}t_{d-1} 
\\ \bangou 
\label{eq4}
p^{(d-2)e}t_{d-2} & = & \frac{1}{2}(c_1^2 - 2c_2)t_d + c_1t_{d-1} + p^{de}t_{d-2} \\
\nonumber
& \vdots &
\end{eqnarray}
Note that $t_d = [X]$, and it is sometimes denoted by $1$ since $[X]$ is the unit
element of the Chow ring of $X$.
By the definition of $t_i$ (Example~3.2.4 and Chapter~18 in \cite{F}), we have
\begin{eqnarray}
\bangou\label{eq5}
t_{d-1} & = & {\rm td}_1(\Omega_X^\vee) = \frac{1}{2}c_1(\Omega_X^\vee)
= - \frac{1}{2}c_1(\Omega_X) = - \frac{1}{2}c_1(\omega_X) = -\frac{1}{2}K_X \\
\bangou\label{eq6}
t_{d-2} & = & {\rm td}_2(\Omega_X^\vee) = 
\frac{1}{12}\left( c_1(\Omega_X^\vee)^2 + c_2(\Omega_X^\vee) \right)
= \frac{1}{12}\left( K_X^2 + c_2(\Omega^\vee_X) \right) .
\end{eqnarray}
Substituting (\ref{eq5}) and (\ref{eq6}) for (\ref{eq3}) and (\ref{eq4}),
we have
\begin{eqnarray*}
c_1(F^e_*{\cal O}_X) & = & \frac{p^{de} - p^{(d-1)e}}{2}K_X \\
c_2(F^e_*{\cal O}_X) & = & 
\frac{3p^{2de} - 6p^{(2d-1)e} + 3p^{2(d-1)e} - 4p^{de} + 6 p^{(d-1)e} - 2p^{(d-2)e}}{24}K_X^2 \\
& & + \frac{p^{de} - p^{(d-2)e}}{12}c_2(\Omega_X^\vee) \\
& \vdots & 
\end{eqnarray*}
We have completed the proof of Theorem~\ref{analogue}.

\section{Some examples}

\bangou
\begin{Example}
\begin{rm}
\begin{enumerate}
\item
This example is due to Han-Monsky~\cite{HM}.
Set 
\[
A = {\Bbb F}_5[[x_1, \ldots, x_4]]/(x_1^4 + \cdots + x_4^4)
\]
 and
$\m = (x_1, \ldots, x_4)A$.
Then, we have
\[
\ell_A(A/\m^{[5^e]}) = \frac{168}{61}5^{3e} 
- \frac{107}{61}3^{e} .
\]
Therefore, in this case, $e_{HK}(\m,A) = \frac{168}{61}$ and $\beta(\m,A) = 0$.
We thus know that there is no hope to extend Theorem~\ref{HMM} under the same assumptions to get a third term in the Hilbert-Kunz function of the form
$\gamma \cdot p^{(d-2)e} + O(p^{(d-3)e})$ in place of $O(p^{(d-2)e})$.
\item
Set 
\bangou
\begin{equation}\label{ring}
A = k[[x_{ij} \mid i = 1, \ldots, m; \ j = 1, \ldots, n]]/I_2(x_{ij}) ,
\end{equation}
where $k$ is a perfect field of characteristic $p > 0$,
and $I_2(x_{ij})$ is the ideal generated by all the $2$ by $2$ minors of 
the generic $m$ by $n$ matrix $(x_{ij})$.

Suppose $m = 2$ and $n = 3$.
Then, K.-i.~Watanabe~\cite{Watanabe} proved 
\[
\ell_A(A/\m^{[p^e]}) = (13p^{4e}-2p^{3e}-p^{2e}-2p^e)/8 .
\]
Therefore, we have $e_{HK}(\m,A) = \frac{13}{8}$ and 
$\beta(\m,A) = -\frac{1}{4} \neq 0$.
\item
Let $(A,\m)$ be a homomorphic image of a regular local ring $S$.
Then, by the singular Riemann-Roch theorem (Chapter~18 and 20 in Fulton~\cite{F}),
we have an isomorphism 
\[
\tau_{\spec(A)/\spec(S)} : \g(A)\subq \longrightarrow \chow{*}(A)\subq
\]
of ${\Bbb Q}$-vector spaces.
For a finitely generated $A$-module $M$, put
\[
\tau_{\spec(A)/\spec(S)}([M]) = \tau_d([M]) + \tau_{d-1}([M]) + \cdots + \tau_0([M]) ,
\]
where $\tau_i([M]) \in \chow{i}(A)\subq$ for $i = 0, 1, \ldots, d$.
Let ${\Bbb F}.$ be a bounded finite $A$-free complex such that each homology module has
 finite length.
Then, by the local Riemann-Roch formula (Example~18.3.12 in \cite{F}),
we have
\[
\sum_{j} (-1)^j \ell_A(H_j({\Bbb F}.\otimes_AM))
= \sum_{i} {\rm ch}({\Bbb F}.) \cap \tau_i([M]) .
\]
Furthermore, assume that $A$ is a Cohen-Macaulay ring of characteristic $p$, where $p$ is a prime number, and
the residue class field of $A$ is perfect.
Let $I$ be a maximal primary ideal of finite projective dimension.
Let ${{\Bbb F}_I}.$ be a finite $A$-free resolution of $A/I$.
Then, if the depth of $M$ is equal to $d$, we have
\begin{eqnarray*}
\ell_A(M/I^{[p^e]}M) & = & \ell_A(F^e_*(M)/IF^e_*(M)) \\
& = & \sum_j(-1)^j \ell_A(H_j(F^e_*(M)\otimes {{\Bbb F}_I}.)) \\
& = & \sum_{i} {\rm ch}({{\Bbb F}_I}.) \cap \tau_iF_*^e([M]) \\
& = & \sum_{i} {\rm ch}({{\Bbb F}_I}.) \cap p^{ie}\tau_i([M]) \\
& = & \sum_{i} \left( {\rm ch}({{\Bbb F}_I}.) \cap \tau_i([M]) \right) p^{ie} 
\end{eqnarray*}
by the equality (\ref{rmn3}).
Therefore, in this case, we have $e_{HK}(I,M) = {\rm ch}({{\Bbb F}_I}.) \cap \tau_d([M])$ and
$\beta(I,M) = {\rm ch}({{\Bbb F}_I}.) \cap \tau_{d-1}([M])$.

One can prove that there exists a maximal primary ideal $I$ 
of finite projective dimension such that 
${\rm ch}({{\Bbb F}_I}.) \cap \tau_{i}([M]) \neq 0$ if and only if
$\tau_{i}([M])$ is not numerically equivalent to $0$ (cf., Theorem~6.4 in \cite{K23}).
We refer the reader to \cite{K23} for the theory of numerical equivalence.

Suppose that $A$ is the ring in (\ref{ring}) as above.
In this case, $\tau_{d-1}([A])$ is not numerically equivalent to $0$ 
if and only if $m \neq n$ (cf., Section~3 in \cite{K11} and Example~7.9 in \cite{K23}).
Therefore, if $m \neq n$, 
then there exists a maximal primary ideal $I$ of finite projective
dimension such that $\beta(I,A) \neq 0$.

On the other hand, if $m = n$, then $A$ is a Gorenstein ring.
By Corollary~\ref{cor}, 
$\beta(I,A) = 0$ for any maximal primary ideal $I$ of $A$.
\end{enumerate}
\end{rm}
\end{Example}

\bangou
\begin{Example}
\begin{rm}
\begin{enumerate}
\item
Set 
\[
A = k[[x_1, x_2, x_3, y_1, y_2, y_3]] \left/
I_2
\left(
\begin{array}{ccc}
x_1 & x_2 & x_3 \\ y_1 & y_2 & y_3
\end{array}
\right) \right. , 
\]
${\frak p} = (x_1, x_2, x_3)A$ and ${\frak q} = (x_1, y_1)A$.
Here, assume that $k$ is a perfect field of characteristic $2$.
Then, applying Hirano's formula\footnote{
Hirano proved the following:
Let $X$ be a $d$-dimensional toric variety over a perfect field 
of characteristic $p > 0$ defined by a fan in $N = {\Bbb Z}^d$.
Let $F : X \rightarrow X$ be the absolute Frobenius map.
Then, for any positive integer $e$,
we have
\[
F^e_*{\cal O}_X = \bigoplus_{0 \leq s_1, \ldots, s_d \leq p^e}
{\cal O}_X
\left(
\frac{1}{p^e}{\rm div}_X(u_1^{s_1}\cdots u_d^{s_d})
\right)
u_1^{s_1/p^e}\cdots u_d^{s_d/p^e} ,
\]
where $\{ u_1, \ldots, u_d \}$ is the dual basis of $N = {\Bbb Z}^d$.
} (Theorem~2 in \cite{Hirano}),
one can prove that
\[
{}^1A \simeq A^{\oplus 10} \oplus {\frak p} \oplus {\frak q}^{\oplus 5} .
\]
Here, remark that $\rank_A{}^1A = p^{\dim A} = 2^4 = 16$.

Then, we have
\[
{\rm cl}({}^1A) = 10 \cdot {\rm cl}(A) + {\rm cl}({\frak p}) + 5 \cdot {\rm cl}({\frak q}) = 4 \cdot {\rm cl}({\frak q}) 
\]
since ${\rm cl}(A) = 0$ and ${\rm cl}({\frak p}) + {\rm cl}({\frak q}) = 0$.

On the other hand, 
it is well known that $\omega_A \simeq {\frak q}$.
By Theorem~\ref{main}, we have
\[
{\rm cl}({}^1A) = \frac{2^4 - 2^3}{2} {\rm cl}(\omega_A) = 4 \cdot {\rm cl}({\frak q}) .
\]
\item
Let $k$ be a perfect field of characteristic $p$,
where $p$ is a prime integer.
Put $X = {\Bbb P}_k^1$.
Let $F : X \rightarrow X$ be the absolute Frobenius map.
Then, we have
$F_*{\cal O}_X \simeq {\cal O}_X \oplus {\cal O}_X(-1)^{\oplus (p-1)}$,
and 
\[
c_1(F_*{\cal O}_X) = c_1(\wedge^p F_*{\cal O}_X) = c_1({\cal O}_X(1-p)) = 1-p .
\]
Here, remark that the natural map
\[
{\rm deg} : {\rm Cl}(X) \rightarrow {\Bbb Z}
\]
is an isomorphism in this case.

On the other hand, it is well known that $\omega_X \simeq {\cal O}_X(-2)$.
Therefore, we have $K_X = -2$.
By Theorem~\ref{analogue}, we have
\[
c_1(F_*{\cal O}_X) = \frac{p - 1}{2} K_X = 1-p .
\]
\end{enumerate}
\end{rm}
\end{Example}

\noindent
{\bf Acknowledgement:} \
The author deeply thanks the referee for his valuable comments.

\noindent
\begin{tabular}{l}
Department of Mathematics \\
Faculty of Science and Technology \\
Meiji University \\
Higashimita 1-1-1, Tama-ku \\
Kawasaki 214-8571, Japan \\
{\tt kurano@math.meiji.ac.jp} \\
{\tt http://www.math.meiji.ac.jp/\~{}kurano}
\end{tabular}


{\small

\begin{thebibliography}{99}

\bibitem{B} {\sc N. Bourbaki},
{\it Alg\`{e}bre commutative, \'{E}l\'{e}ments de Math\'{e}matique, Chap.\ 1--7},
Hermann Paris, 1961--1965.


\bibitem{F} {\sc W.~Fulton},
{\it Intersection Theory, 2nd Edition},
Springer-Verlag, Berlin, New York, 1997.


\bibitem{HM} {\sc C. Han and P. Monsky},
{\it Some surprising Hilbert-Kunz functions},
Math.\ Zeit.\ {\bf 214} (1993), 119--135.


\bibitem{Hirano} {\sc N. Hirano},
{\it Frobenius direct images and Hilbert-Kunz functions (in Japanese)},
Abstract of the second Kyoto COE seminar in Kinosaki,
Feburary, 2005, to appear.


\bibitem{HMM} {\sc C. Huneke, M. A. McDermott and P. Monsky},
{\it Hilbert-Kunz functions for normal graded rings},
Math.\ Res.\ Letters {\bf 11} (2004), 539--546.





\bibitem{K9} {\sc K. Kurano},
{\it An approach to the characteristic free Dutta multiplicities},
J. Math.\ Soc.\ Japan, {\bf 45} (1993), 369--390.

\bibitem{K11} {\sc K. Kurano},
{\it A remark on the Riemann-Roch formula for affine schemes associated with 
Noetherian local rings}, T\^{o}hoku Math.\ J. {\bf 48} (1996), 121--138. 




\bibitem{K23} {\sc K. Kurano},
{\it Numerical equivalence defined on Chow groups of Noetherian local rings},
Invent.\ Math.\ {\bf 157} (2004), 575--619.








\bibitem{Mo} {\sc P. Monsky},
{\it The Hilbert-Kunz function},
Math.\ Ann.\ {\bf 263} (1983), 43--49.
















\bibitem{Watanabe} {\sc K.-i. Watanabe},
a private discussion.
\end{thebibliography}
}
\end{document}